\documentclass[11pt]{article}
\usepackage{amsfonts,amsmath}
\usepackage{amsmath, amsthm, upref}
\usepackage{graphicx}

\newtheorem{theorem}{Theorem}
\newtheorem{lemma}[theorem]{Lemma}
\newtheorem{lem}[theorem]{Lemma}
\newtheorem{prop}[theorem]{Proposition}

{\bf}{\it}
{\bf}{\it}
\newtheorem{corollary}{Corollary}

\newcommand{\mcal}{\mathcal}
\newcommand{\e}{E}
\newcommand{\rr}{\mathbb{R}}
\newcommand{\comment}[1]{}

\newcommand{\abs}[1]{\left\lvert #1\right\rvert}

\newcommand{\gbar}{\bar{g}}
\newcommand{\iprod}[1]{\left\langle #1 \right\rangle}

\newcommand{\tp}{P}

\newcommand{\eq}{\begin{equation}}
\newcommand{\en}{\end{equation}}

\begin{document}

\theoremstyle{plain}
\newtheorem{thm}{Theorem}

\graphicspath{/Users/soumikpal/Documents/Academics/Research/Protter}

\title{Analysis of continuous strict local martingales via h-transforms}
\author{Soumik Pal\thanks{Supported in part by NSF grant DMS-0306194;  Departments of Mathematics -- C-547 Padelford Hall, University of Washington, Seattle 98195 USA}\; and Philip Protter
\thanks{Supported in part by NSF grant DMS-0202958 and NSA grant MDA-904-03-1-0092; ORIE -- 219
Rhodes Hall,
Cornell University, Ithaca, NY 14853-3801 USA} }

\date{\today}
\maketitle

\begin{abstract}
We study strict local martingales via $h$-transforms, a method which
first appeared in Delbaen-Schachermayer. We show that strict local
martingales arise whenever there is a consistent family of change of
measures where the two measures are not equivalent to one another.
Several old and new strict local martingales are identified. We
treat examples of diffusions with various boundary behavior,
size-bias sampling of diffusion paths, and non-colliding diffusions.
A multidimensional generalization to conformal strict local
martingales is achieved through Kelvin transform. As curious
examples of non-standard behavior, we show by various examples that
strict local martingales do not behave uniformly when the function
$(x-K)^+$ is applied to them. Implications to the recent literature
on financial bubbles are discussed.
\end{abstract}

\section{Introduction}
Local martingales which are not martingales (known as ``strict''
local martingales) arise naturally in the Doob-Meyer decomposition
and in multiplicative functional decompositions, as well as in
stochastic integration theory. They are nevertheless often
considered to be anomalies, processes that need to be maneuvered by
localization. Hence, studies focussed purely on strict local
martingales are rare. One notable exception is the article by
Elworthy, Li, \& Yor \cite{ELY} who study their properties in depth.
On the other hand, applications of strict local martingales are
common. See the articles by Bentata \& Yor \cite{bentatayor}, Biane
\& Yor \cite{biyor88}, Cox \& Hobson \cite{CH}, Fernholz \& Karatzas
\cite{ferkar05}, and the very recent book-length preprint of
Profeta, Roynette, and Yor~\cite{PRY}.

Our goal in this paper is to demonstrate that strict local
martingales capture a fundamental probabilistic phenomenon.
Suppose there is a pair of probability measures on a filtered probability space, where
one strictly dominates the other in the sense that their null sets
are not the same. Then one can establish a correspondence between pairs of continuous processes such that one is a martingale under one measure and the other is a strict local martingale under the changed measure. This has been made precise in Propositions \ref{mainthm} and \ref{htrans}.
The continuity of the time parameter is crucial since there are no strict local martingales in discrete time.

The above idea is not new and has appeared in the literature under various disguises. For example, for positive local martingales, such a phenomenon was identified by Delbaen and Schachermayer in \cite{ds95}. 
We, however, do a systematic exploration to identify several strict local martingales as an application of the previous result, in diverse topics such as diffusions
conditioned to exit through a subset of the boundary of a domain,
size-biased sampling of diffusion paths, and non-colliding
diffusions such as Dyson's Brownian motion from the Random Matrix
Theory. In addition, we show how usual stochastic calculus with respect to strict local martingales can be made more friendly by incorporating the change of measure. Our results in this direction are related to recent work by Madan and Yor \cite{MY}.

We also prove a multidimensional analogue of our results where a
strict local martingale in one-dimension is replaced by a conformal
local martingale in three or more dimensions where at least one
coordinate process is strict. The analysis exploits the Kelvin
transform from classical potential theory.

The name \textit{h-transforms} in the title refers to the classical technique (originally due to Doob), where the law of a Markov process is changed by using a nonnegative harmonic function. Aided by Girsanov's theorem, this usually leads to a description of the process conditioned to stay away from the zero-set of the harmonic function. As our examples will demonstrate, this is a fertile set-up for a plethora of strict local martingales. For more details about classical $h$-transforms, please see Revuz and Yor \cite[p.~357]{ry99}.

The rest of the paper analyzes the effect of convex functions applied to strict local martingales. Although these are always local submartingales, the expectations of the processes can be of arbitrary shapes. Similar examples were also identified in \cite{ELY}. In particular, this has some consequences in mathematical finance.
A natural source for local martingales in mathematical finance is the condition of No Free Lunch with Vanishing Risk (see Delbaen and Schachermayer \cite{DS}). Roughly, it states that in a financial market the no arbitrage condition is
equivalent (in the case of continuous paths) to the existence of an
(equivalent) ``risk neutral'' probability measure $Q$ which turns
the price process into either a martingale or a strict local
martingale. The implications of our results can be readily
understood if we assume that the risk neutral measure produces a
(one-dimensional) price process $X=(X_t)_{t\geq 0}$ that is a strict
local martingale. In that case, for example, the process $Y_t=(X_t-K)^+$
\emph{need not be a submartingale}, and the function $t\mapsto
E\{(X_t-K)^+\}$ need no longer be increasing, contradicting the
usual wisdom in the theory. The original purpose of this paper was
to understand this phenomenon better, motivated in particular by the
role local martingales play in financial bubbles (cf~\cite{JPS1}
and~\cite{JPS2}).

\section{A method to generate strict local martingales}

Let $\left(  \Omega, \left\{ \mathcal{F}_t\right\}_{t\ge 0} \right)$
be a filtered sample space on which two probability measures $P$ and
$Q$ are defined. We assume that $P$ is locally strictly dominated by
$Q$, in the sense that $P$ is absolutely continuous with respect to
$Q$ $(P \ll Q)$ on every $\mcal{F}_t$ (see \cite[Chap.~VIII]{ry99}). Let $h_t$ denote the Radon-Nikod\'ym process ${dP}/{dQ}$ on the $\sigma$-algebra ${\mathcal{F}_t}$. Assume that \eq\label{locstrdom} Q\left( \tau_0 < \infty \right) > 0, \en where
$\tau_0 = \inf\left\{ s > 0, \; h_s =0 \right\}$. We claim the
following result.

\begin{prop}\label{mainthm}
Assume that $h$ is a continuous process $Q$-almost surely. Let
$\{f_t, \; t \ge 0\}$ be a continuous $Q$-martingale adapted to the
filtration $\{ \mathcal{F}_t\}$. Suppose either of the two
conditions hold:
\begin{enumerate}
\item[(i)] $f$ is uniformly integrable, $E^Q(f_0)\neq 0$, and $Q\left( \tau_0 < \infty  \right)=1$, or,
\item[(ii)] $f$ is nonnegative, $Q(\sigma_0 > \tau_0)> 0$, where $\sigma_0=\inf\{ s\ge 0: f_s=0 \}$.
\end{enumerate}
Then the process $N_t := {f_t}/{h_t}$ is a strict local martingale under $P$.
However, if $f$ is nonnegative and $Q(\sigma_0 > \tau_0)=0$, the
process $N_t$ is a martingale.
\end{prop}

\begin{proof}[Proof of Proposition \ref{mainthm}] Let us first show
that $N$ is a local martingale. Consider the sequence of stopping
times
\[
\tau_k := \inf \left\{ s\ge 0: \; h_s \le 1/k   \right\},\quad k=1,2,\ldots.
\]
Then, it follows that $P\left( \lim_{k\rightarrow \infty} \tau_k
= \infty  \right)=1$ (see \cite[p.~326]{ry99}). 
Thus, for any bounded stopping time $\tau$, we get by the change of measure
formula
\[
E^P\left(  N_{\tau \wedge \tau_k}  \right)= E^Q\left(
h_{\tau\wedge \tau_k} \frac{1}{h_{\tau\wedge
\tau_k}}f_{\tau\wedge \tau_k}\right)= E^Q\left( f_{\tau\wedge\tau_k} \right)= E^Q\left( f_{0} \right).
\]
The final equality is due to the fact that $f$ is a $Q$-martingale and $\tau\wedge \tau_k$ is a bounded stopping time.

Since $ N_{\tau \wedge \tau_k}$ has the same expectation for all
bounded stopping times $\tau$, it follows that $N_{\cdot\wedge
\tau_k}$ is a martingale. The local martingale property now
follows.

To show now that it is not a martingale, we compute the expectation
of $N_t$. Note that $h$ is a nonnegative $Q$-martingale. Thus, again applying the
change of measure formula, we get \eq\label{ptoq} E^P\left( N_t
\right) = E^Q\left(  f_t1_{\left\{  \tau_0 > t \right\}}  \right),
\en where $\tau_0$ is the hitting time of zero for $h$.

Now suppose (i) $\{f_t\}$ is a uniformly integrable martingale and
$Q\left( \tau_0 < \infty  \right)=1$. By uniformly integrability we
see from \eqref{ptoq} that
\[
\lim_{t\rightarrow \infty} E^P(N_t)=0
\]
which shows that $N$ is not a martingale, since $E^P(N_0)=E^Q(f_0)$ is
assumed to be non-zero.

Finally suppose (ii) $Q(\sigma_0 > \tau_0) > 0$ holds. Since zero is
an absorbing state for the nonnegative martingale $f$, by
\eqref{locstrdom} there is a time $t >0$ such that $Q(\{ f_t > 0\}
\cap \{ \tau_0 \le t \}) > 0$. For that particular $t$, we get from
\eqref{ptoq} that
\[
E^P\left( N_t \right) < E^Q(f_t) = E^Q(f_0) = E^P(N_0).
\]
This again proves that the expectation of $N_t$ is not a constant.
Hence it cannot be a martingale.

The final assertion follows from \eqref{ptoq} by noting that $N$, by
virtue of being a nonnegative local martingale, is a supermartingale
which has a constant expectation. Thus it must be a martingale.
\end{proof}

It is clear that some condition is necessary for the last theorem to
hold as can be seen by taking $f_t\equiv h_t$ which results in a
true martingale $N_t\equiv 1$.  

Several examples of old and new strict local martingales follow from
the previous result by suitably choosing a change of measure and the
process $f$. We describe some classes of examples below.

\subsection{Diffusions with different boundary behaviors}

One of the earliest uses of change of measures was to condition a
diffusion to exhibit a particular boundary behavior. This typically
involves a change of measure that is locally strictly dominated. We
give below two examples to show how local martingales arise from
such a set-up.
\bigskip

\noindent\textbf{Example 1.} The $3$-dimensional Bessel process,
BES($3$), (see Karatzas \& Shreve \cite[page 158]{KS} for the
details) is the (strong) solution of the stochastic differential
equation:
\begin{equation}\label{bessde}
dX_t = \frac{1}{X_t}dt + d\beta_t, \qquad X_0=x_0 \ge 0,
\end{equation}
where $\beta$ is an one-dimensional standard Brownian motion. It is
also the law of the Euclidean norm of a three-dimensional Brownian
motion.  The reciprocal of this process, known as the \textit{inverse Bessel process},
serves as a prototypical example of a local martingale which is not
a martingale. 

This is an immediate example of Proposition \ref{mainthm} where we consider the canonical sample space for a Brownian motion, define $Q$ to be law of the Brownian motion, starting at $x_0$, and killed at zero, and let $h_t$ be the coordinate map $\omega_t/x_0$. Finally we take $f_t\equiv 1$. Details on why $Q$ is the law of the Bessel process are well-known and can be found in \cite[p.~252]{ry99}.

Several other examples of the same spirit can be derived. In
particular, for any Bessel process $X$ of dimension $\delta > 2$, it
is known that $X^{2-\delta}$ is a strict local martingale. This
can be proved similarly as above using the martingale
$Y^{\delta-2}$, where $Y$ is a Bessel process of dimension
$(4-\delta)$ (which can be negative) absorbed at zero. See the
article \cite{bentatayor} for the details.

\bigskip

\noindent\textbf{Example 2.} The previous example can be easily
extended to a multidimensional form. Let $D$ be an open, connected, bounded
domain in $\rr^n$ ($n \ge 2$) where every point on the boundary is
regular (in the sense of \cite[p. 245]{KS}). Consider an
$n$-dimensional Brownian motion $X$ starting from a point $x_0 \in
D$ getting absorbed upon hitting the boundary of $D$, say $\partial
D$. Let $Q$ denote the law of the process $\{ X_{t\wedge \tau_D}, \;
t \ge 0 \}$. Since $D$ is bounded, $Q(\tau_D < \infty)=1$. For any bounded measurable nonnegative function $u$ on
$\partial D$, one can construct the following function
\eq\label{harext} f(x_0)= E^Q\left(  u\left( X_{\tau_D}\right) \mid
X_0=x_0   \right). \en 
 In fact, $f$ is the solution
of the Dirichlet problem on $D$ with boundary data $u$.
By the Markov property, it follows that
$f(X_{t\wedge\tau_D})$ is a martingale.

Let $B_1$ be a connected (nontrivial) proper subset of the boundary
$\partial D$, and let the function $u$ be one on  $B_1$ and zero
elsewhere, i.e, $u(x)=1_{(x\in B_1)}$. In that case, the resulting
harmonic function in $D$ is given by
\[
v(x_0)=Q\left( X_{\tau_D} \in B_1 \mid X_0=x_0  \right).
\]
Then $v(x_0)>0$, and by the tower property, $h(X_{t\wedge \tau_D})= v(X_{t\wedge
\tau_D})/v(x_0)$ is a nonnegative martingale starting from one. Let
$P$ denote the change of measure using the Radon-Nikod\'ym derivative process $h(X_{t\wedge \tau_D})$. This gives us a process law which can be interpreted as
\textit{Brownian motion, starting from $x_0$, conditioned to exit through $B_1$}.

As a corollary of Proposition \ref{mainthm}, we get the following result.

\begin{prop}
Let $u$ be any nonnegative function on $\partial D$ such that $u >
0$ on some subset of $\partial D \backslash B_1$ of positive
Lebesgue measure. Let $f$ be the harmonic extension of $u$ given by
 \eqref{harext}. Then the process
\[
N_t = \frac{f(X_{t\wedge \tau_D})}{v(X_{t\wedge \tau_D})}
\]
is a strict local martingale under $P$.
\end{prop}

\begin{proof}
Follows from Proposition \ref{mainthm} by verifying condition
(ii) since $u(X_{\tau_D})$ need not always be zero when $X_{\tau_D} \notin B_1$.
\end{proof}

\subsection{Size-biased sampling of diffusion paths}

Another class of interesting examples follow in the case of the
size-biased change of measure. Size-biased sampling has been often
discussed in connection with discrete distributions, see, for
example, the article by Pitman \cite{pitman96} and the references
within. It usually involves a finite or countable collection of
numbers $\{ p_1, p_2, p_3, \ldots \}$ such that each $p_i \ge 0$ and
$\sum_i p_i=1$. A size-biased sample from this collection refers to
a sampling procedure where the sample $\tilde p$ has the
distribution
\[
P\left(  \tilde p = p_i \right)=p_i, \quad \text{for all}\; i=1,2,\ldots.
\]
One can now remove this chosen sample from the collection,
renormalize it, and repeat the procedure. This is closely connected
to urn schemes where each $p_i$ refers to the proportion of balls of
a color $i$ that is in an urn. If one randomly selects a ball, the
color of the chosen ball has the size-biased distribution.

One can similarly develop a concept of size-biased sampling of
diffusion paths as described below. Consider $n$ non-negative
diffusions $\{ X_1, \ldots, X_n\}$ running in time. Fix a time $t$,
and look at the paths of the diffusions during the time interval
$[0,t]$. Denote these random continuous paths by $X_1[0,t],
X_2[0,t], \ldots, X_n[0,t]$. Sample one of these random paths with 
probability proportional to the terminal value $X_i(t)$. That is,
the sampled path has the law
\[
Y[0,t] = X_i[0,t], \quad \text{with probability}\; \frac{X_i(t)}{\sum_{j=1}^n X_j(t)}.
\]
How can one describe the law of $Y$ ? In general the law of $Y$
might not be consistently defined as time varies. Nevertheless there
are cases where it makes sense. We show below the example when each
$X_i$ is a Bessel square process of dimension zero (BESQ$^0$). This
is the strong solution of the SDE
\[
Z_t = z + 2\int_0^t \sqrt{Z_s} d\beta_s, \qquad z >0,
\]
where $\beta$ is a one-dimensional Brownian motion. BESQ$^0$ is a
nonnegative martingale also known as the Feller branching diffusion,
in the sense that it represents the total surviving population of a
critical Galton-Watson branching process. Indeed, our treatment here
of size-biased transforms of BESQ processes is inspired by
size-biased transforms of Galton-Watson trees (see the work by
Lyons-Pemantle-Peres \cite{LPP} and references to the prior
literature referred to there).

The interpretation of such dynamic size-biased sampling when every
$X_i$ is a BESQ$^0$ is straightforward. When each $X_i$ has the law
of Feller's branching diffusion, they represent a surviving
population from $n$ critical branching processes. We can do
size-biased sampling at different time points from these
populations. The construction below describes the joint law of these
samples.

Consider the canonical sample space for multidimensional diffusions,
i.e., the $n$-dimensional continuous path space $C^n[0,\infty)$,
coupled with the usual right-continuous and complete (with respect to the diffusion we describe below) filtration
generated by the coordinate maps. Let $\omega_t=(\omega_t(1),
 \ldots, \omega_t(n))$ denote a sample path. Also, as
before, we take $\omega[0,t]$ to denote the path during
time-interval $[0,t]$. Let $Q$ be the joint law of $n$ independent
BESQ$^0$ processes starting from positive points $(z(1), \ldots,
z(n))$. To keep matters simple, we assume all the $z(i)'s$ are the
same and equal to $z$. This induces the following exchangeability
property.

Let $T_i:C^n[0,\infty)\rightarrow C^n[0,\infty)$ be a permutation operator on
the sample space such that $T_i\omega(1)= \omega(i)$ and
\eq\label{whatistnu}
\begin{split}
T_i\omega(j+1)=\omega(j), \quad \text{if}\; j < i, \quad \text{and},
\quad T_i\omega(j)=\omega(j), \qquad \text{if}\; j > i.
\end{split}
\en Thus, $T_i$ puts the $i$th coordinate as the first, and shifts
the others appropriately. By our assumed exchangeability under $Q$,
each $T_i\omega$ has also the same law $Q$. To define the
size-biased change of measure we need the following lemma.

\begin{lemma}\label{ratiomart}
For any $i$, the process
\[
M_t(i) = \frac{\omega_t(i)}{\omega_t(1) + \ldots + \omega_t(n)},
\qquad t \ge 0,
\]
is a martingale under $Q$.
\end{lemma}

\begin{proof}
Without loss of generality take $i=1$. Let $Z(1), Z(2), \ldots,
Z(n)$ be independent BESQ$^0$. The sum $\zeta=Z(1) + \ldots +Z(n)$
is another BESQ$^0$ process. Thus by It\^o's rule we get
\[
\begin{split}
d\left(  {Z_t(1)}/{\zeta_t} \right)&= Z_t(1) d\left(
\frac{1}{\zeta_t}
\right) + \frac{1}{\zeta_t}dZ_t(1) + d\iprod{Z_t(1), \zeta_t^{-1}}\\
&= \text{local martingale}\; + \frac{Z_t(1)}{\zeta_t^3}4\zeta_tdt - \frac{4Z_t(1)}{\zeta_t^2}dt.
\end{split}
\]
This proves that the ratio process is a local martingale. But since
it is bounded, it must be a martingale.
\end{proof}

Define the size-biased sampling law $\tp$ on $(C^n[0,\infty),
\{\mcal{F}_t\})$ by \eq\label{whatistpnu} \tp(A)= E^Q\left[
\sum_{i=1}^n  \frac{\omega_t(i)}{\omega_t(1)+ \cdots+\omega_t(n)}
1_{\left\{ T_i\omega[0,t] \in A \right\}} \right], \quad \text{for
all}\; A\in \mcal{F}_t. \en Note that, by Lemma \ref{ratiomart},
$\tp$ defines a consistent probability measure on the filtration
$\{\mcal{F}_t\}$. We show below that $P$ is strictly
locally dominated by $Q$ and compute the Radon-Nikod\'ym derivative.

Recall that, under $Q$, each $T_i\omega$ has the same
law $Q$. Thus, we can simplify expression \eqref{whatistpnu} to
write
\[
\tp(A)= n E^Q\left[
\frac{\omega_t(1)}{\omega_t(1)+\cdots+\omega_t(n)} 1_{\left\{
\omega[0,t] \in A \right\}} \right], \quad \text{for all}\; A\in
\mcal{F}_t.
\]
This proves that $P \ll Q$ and the Radon-Nikod\'ym derivative is given by
\[
h_t= \frac{n\omega_t(1)}{\omega_t(1)+\cdots+\omega_t(n)}.
\]
Since under the BESQ$^0$ law, every coordinate can hit zero and get
absorbed, the above relation is a locally strict domination and
hence leads to examples of strict local martingales.

It can be verified, using Girsanov Theorem, that under $P$, the coordinate processes remain independent. The first coordinate $\omega(1)$ has law BESQ of dimension $4$, and the rest of the coordinates are still BESQ$^0$. These correspond to the spine decomposition of critical GW trees, conditioned on surviving forever, as in \cite{LPP}.

As an immediate corollary of Proposition \ref{mainthm} we get the following result.

\begin{prop} Let $(Z(1), \ldots, Z(n))$ be continuous processes
whose law is the size-biased sampled BESQ$^0$ law $P$ described
above. Let $\zeta_t=Z_t(1)+\ldots + Z_t(n)$ be the total sum
process. Then the processes
\[
N_t= \frac{\zeta_t^2}{ Z_t(1)},\quad U_t =   \frac{Z_t(2)\zeta_t} {
Z_t(1)},\quad V_t= \frac{\zeta_t}{Z_t(1)}\prod_{i=2}^n Z_t(i),\quad
t \ge 0,
\]
are all strict local martingales. However, the process
\[
M_t = \zeta_t\prod_{i=2}^n Z_t(i),\quad t \ge 0,
\]
is a true martingale.
\end{prop}

\begin{proof}
Let $Z(1), \ldots, Z(n)$ be iid BESQ$^0$ processes starting from $z
>0$. Then $\zeta$, $Z(2)$, and $\prod_{i=2}^n Z_t(i)$ are all true
martingales which can remain positive when $Z_t(1)=0$. The result
now follows from Proposition \ref{mainthm}, condition (ii).

On the other had, the process $\prod_{i=1}^n Z_t(i)$ is a martingale
that always hits zero before $Z_t(1)$. Thus, $M$ is a true
martingale.
\end{proof}

\subsection{Non-colliding diffusions} Our third class of examples
are cases when the change of measure leads to non-intersecting paths
of several linear diffusions. Probably the most important example of
this class is Dyson's Brownian motion which is a solution of the
following $n$-dimensional SDE: \eq\label{dysonsde} d\lambda_t(i) =
\sum_{j\neq i} \frac{2}{\lambda_t(i) - \lambda_t(j)} dt +  dB_t(i),
\quad t\ge 0, \quad i=1,\ldots,n. \en Here $(B(1), \ldots, B(n))$ is
an $n$-dimensional Brownian motion. It appears in the context of
Random Matrix Theory. Please see the survey article by K\"onig
\cite{konigsurvey} and the original paper by Dyson \cite{dysonpaper}
for the details (including the definition of the Gaussian Unitary
Ensemble) and the proofs.

\begin{thm}
For any $i=1,2,\ldots,n$, and $i< j\le n$, let $\{M_t(i,i), \; t \ge 0
\}$, $\{ M^R_t(i,j), \; t \ge 0 \}$, and $\{ M^I_t(i,j), \; t\ge 0
\}$ be independent real standard Brownian motions, starting at zero.
The Hermitian random matrix $M_t=(M_t(i,j), \; 1\le i,j\le n)$, with
\[
 M_t(i,j)=M_t^R(i,j) + i M_t^I(i,j),\quad i < j,
\]
has the distribution of the Gaussian Unitary Ensemble at time $t=1$.
Then the process $(\lambda_t, \; t \ge 0)$ of $n$ eigenvalues of
$M_t$ satisfies SDE \eqref{dysonsde}. It can interpreted as a
conditional Brownian motion in $\rr^n$, starting at zero,
conditioned to have
\[
\lambda_t(1) < \lambda_t(2) < \cdots < \lambda_t(n), \qquad \text{for all}\quad t > 0.
\]
\end{thm}

What is interesting is that the process in \eqref{dysonsde} can be
obtained as a strict local domination relation from the
$n$-dimensional Wiener measure using the harmonic function
\[
\Delta_n(x)= \prod_{1 \le i < j \le n} (x(j) - x(i)), \qquad x=(x(1), \ldots, x(n)),
\]
which is the well-known Vandermonde determinant. That the function
$\Delta_n$ is harmonic can be found in \cite[p. 433]{konigsurvey}
where it is shown that $\Delta_n(W_t)$ is a martingale when $W$ is
an $n$-dimensional Brownian motion and that the law of the process
in \eqref{dysonsde} (say $P$) can be obtained from the
$n$-dimensional Wiener measure $Q$, by using $h_t=\Delta_n(W_t)$ as
the Radon-Nikod\'ym derivative. Since $\Delta_n(W)$ is zero whenever
any two Brownian coordinates $W(i), W(j)$ are equal (``collide''),
we are in the scenario of Proposition \ref{mainthm}. We prove the
following result.

\begin{prop}
Consider Dyson's Brownian motion $(\lambda_t(1), \ldots,
\lambda_t(n))$ in \eqref{dysonsde}, and for some $m < n$ consider
the process
\[
\Delta_m(t) = \prod_{1 \le i < j \le m} (\lambda_t(j) - \lambda_t(i)), \qquad t \ge 0.
\]
Then the process $N_t= \Delta_m(t) / \Delta_n(t)$ is a strict local
martingale. As a consequence if we consider the Vandermonde
matrix-valued process
\[
A_t=(\lambda_t^{j-1}(i), \; i,j=1,2,\ldots,n), \qquad t \ge 0,
\]
Then every process $A^{-1}_t(n,i)$ is a strict local martingale for $i=1,2,\ldots,n$.

\end{prop}

\begin{proof} Note that, if $W$ is an $n$-dimensional Brownian motion
then $\Delta_m(W_t)$, for any $m\le n$, is a true martingale. Since
$m < n$, it is possible to have $\Delta_m(W_t)$ to be positive when
$\Delta_n(W_t)=0$. The result that $N$ is a strict local martingale
now follows from Proposition \ref{mainthm}.

For the second assertion we use adjugate (or, classical adjoint) relation for the inverse. Let $\text{det}(B)$ for a square matrix $B$ refer to its
determinant. Then
\[
A_t^{-1}(n,i)= \frac{(-1)^{i+n}}{\text{det}(A_t)} \text{det}
\left( \hat{A}_t(i,n)\right), \quad i=1,2,\ldots,n,
\]
where $\hat{A}_t(i,n)$ is the matrix obtained from $A_t$ by removing
the $i$th row and the $n$th column.

Now, $A_t$ is the Vandermonde matrix, so its determinant is equal to
$\Delta_n(t)$. If we remove the $i$th row and the $n$th column from
$A_t$ we get an $(n-1)\times (n-1)$ order Vandermonde matrix of all
the $\lambda_j$'s except the $i$th. Its determinant is again a
Vandermonde determinant. Each of the ratios
\[
\frac{\text{det}\left( \hat{A}_t(i,n)\right)}{\text{det}(A_t)} ,
\]
is a strict local martingale by our earlier argument (for $m=n-1$).
This completes the proof.
\end{proof}

As a final area of applications of $h$-transforms which are locally
strictly dominated, let us mention the theory of measure-valued
processes, in particular, the Dawson-Watanabe superprocesses. A
well-known example is conditioning a superprocess to survive
forever, which can be done by changing the law of a superprocess by
using the total mass process (which is a martingale) as the
Radon-Nikod\'ym derivative. Please see the seminal article in this
direction by Evans \& Perkins \cite{evansperkins}, the book by
Etheridge \cite{etheridge}, and a follow-up article on similar other
$h$-transforms by Overbeck \cite{overbeck}.

\section{A converse to the previous result} As a converse to
Proposition \ref{mainthm} it turns out that all strict local
martingales which remain strictly positive throughout can be
obtained as the reciprocal of a martingale under an $h$-transform.
This was essentially proved by Delbaen and Schachermayer \cite{ds95}
in 1995 in their analysis of arbitrage possibilities in Bessel
processes. We replicate their theorem below. The construction is
related to the F\"ollmer measure of a positive supermartingale
\cite{foll72}.

Before we state the result we need a technique which adds an extra
absorbing point, $\infty$, to the state space $\rr^+$,
originally inspired by the work of P. A. Meyer \cite{meyer72}. We
follow closely the notation used in \cite{ds95}. The space of
trajectories is the space $C_{\infty}[0,T]$ or
$C_{\infty}[0,\infty)$ of continuous paths $\omega$ defined on the
time interval $[0,T]$ or $[0,\infty)$ with values in $[0,\infty]$
with the extra property that if $\omega(t)=\infty$, then
$\omega(s)=\infty$ for all $s > t$. The topology endowed is the one
associated with local uniform convergence. The coordinate process is
denoted by $X$, i.e., $X(t)=\omega(t)$.

\begin{thm}[Delbaen and Schachermayer, Theorem 4 in \cite{ds95} ]\label{dsthm}
If $R$ is a measure on $C[0,1]$ such that $X$ is a strictly positive
strict local martingale, then
\begin{enumerate}
\item[(i)] there is a probability measure $R^*$ on $C_{\infty}[0,1]$ such
that $M=1/X$ is an $R^*$ martingale.
\item[(ii)] We may choose $R^*$ in such a way that the measure $R$ is absolutely
continuous with respect to $R^*$ and its Radon-Nikod\'ym derivative is given by $dR=M_1dR^*$.
\end{enumerate}
\end{thm}

The following result is a corollary.

\begin{prop}\label{htrans} Let $R$ be a probability measure on $C[0,\infty)$
under which the coordinate process $X$ is a positive strict local
martingale starting from one. Then there exists a probability
measure $Q$ on the canonical space such that $X$ is a nonnegative
martingale under $Q$ and the following holds:
\begin{enumerate}
\item[(i)] The probability measure defined by
\begin{equation}\label{abscont}
P(A):= E^Q\left( X_{t} 1_A  \right), \quad \forall\; A\in \mathcal{F}_t,\; t \ge 0,
\end{equation}
is the law of the process $\{1/X_t, \; t \ge 0\}$ under $R$.
\item[(ii)] $X$ is a strict local martingale if and only if $Q(\tau_0 < \infty) > 0$,
where $\tau_0= \inf\{t\ge 0: \; X_t=0\}$ is the first hitting time of zero.
\end{enumerate}
\end{prop}

\begin{proof} Let us first construct a probability measure $Q_1$ on $C[0,1]$ such that \eqref{abscont} holds for $0\le t \le 1$. To wit, consider the restriction of $R$ on $C[0,1]$. The coordinate process, under the restricted measure, remains a positive strict local martingale. Thus, there is an $R^*$ such that Theorem \ref{dsthm} holds. 
Define $Q_1$ to be law of the process $\{ M_t, \; 0\le t \le 1\}$ under $R^*$. Then, we claim that the finite dimensional distributions of the coordinate process $X$ under $Q_1$ is the same as the law of $1/X$, under $R$. To see this, consider $m \in \mathbb{N}$ many points $0\le t_1 < t_2 < \ldots t_{m-1} < t_m\le 1$. For any choice of nonnegative measurable function $g:\rr^m\rightarrow \rr^+$, we get
\eq\label{abscont1}
\begin{split}
E^R\left[ g\left( X_{t_1}^{-1}, X_{t_2}^{-1}, \ldots, X_{t_m}^{-1}  \right) \right]&= E^{R^*}\left[  M_1g\left( M_{t_1}, M_{t_2}, \ldots, M_{t_m}  \right)   \right]\\
&= E^{Q_1} \left[  X_1g\left( X_{t_1}, X_{t_2}, \ldots, X_{t_m}  \right)   \right].
\end{split}
\en
The first equality above is due to part (ii) of Theorem \ref{dsthm}. This proves \eqref{abscont} for $0\le t \le 1$. Since a measure on $C[0,1]$ is determined by its finite dimensional distributions, we are done. 

For any other time $t$, consider the scaled process $\{ X_{tu}, \; 0\le u \le 1\}$ under $R$. It remains a strictly positive strict local martingale, now running between time zero and one. As in the previous paragraph one can construct a measure, $\mu^t_1$, such that, for any choice of $ 0\le t_1 < t_2 <\ldots < t_m \le t$, we have   
\eq\label{abscont2}
E^R\left[ g\left( X_{t_1}^{-1}, X_{t_2}^{-1}, \ldots, X_{t_m}^{-1}  \right) \right] = E^{\mu^t_1}\left[ X_1g\left( X_{t_1/t}, X_{t_2/t}, \ldots, X_{t_m/t}  \right)    \right].
\en
Let $Q_t$ denote the law on $C[0,t]$, such that the law of the scaled coordinate process $\{ X_{tu}, \; 0\le u \le 1\}$ under $Q_t$ is $\mu^t_1$.  Once we demonstrate that this tower of probability measures $\{ Q_t, \; t\ge 0 \}$ is consistent, it follows from standard arguments that they induce a probability measure $Q$ on the entire space $C[0,\infty)$ such that $Q$ restricted $C[0,t]$ is $Q_t$.

To see consistency, let $s < u$ be two positive numbers, and we show that $Q_{u}$ restricted to $\mathcal{F}_{s}$ is $Q_s$. As before, choose $m \in \mathbb{N}$, and consider a nonnegative measurable function $g:\rr^m \rightarrow \rr^+$. Choose
\[
0\le t_1 < t_2 < \ldots t_{m-1} < t_m\le s < u.
\]
Since $X$ is a martingale under $Q_u$ (and $Q_s$), we have 
\[
\begin{split}
E^{Q_{u}}\left[ X_{t_m} g\left(  X_{t_1}, \ldots, X_{t_m} \right) \right]&= E^{Q_{u}}\left[ X_u g\left(  X_{t_1}, \ldots, X_{t_m} \right) \right]=\\
 E^{\mu^u_1}\left[ X_1 g\left(  X_{t_1/u}, \ldots, X_{t_m/u} \right) \right]&= E^{R}\left[ g\left(  X^{-1}_{t_1}, \ldots, X^{-1}_{t_m} \right)    \right], \quad \text{by \eqref{abscont2}},\\
&= E^{\mu^s_1}\left[  X_1 g\left(  X_{t_1/s}, \ldots, X_{t_m/s} \right) \right], \quad \text{again by \eqref{abscont2}},\\
&= E^{Q_{s}}\left[ X_s g\left(  X_{t_1}, \ldots, X_{t_m} \right) \right]=E^{Q_{s}}\left[ X_{t_m} g\left(  X_{t_1}, \ldots, X_{t_m} \right) \right].
\end{split}
\]
Since $X_{t_m} g\left(  X_{t_1}, \ldots, X_{t_m} \right)$ is another function of finite dimensional distributions, the above identity shows the following. For any function $h:(\rr^{+})^m \rightarrow \rr^+$ such that 
\eq\label{coffun}
h(x)=0, \quad \text{if} \quad x_m=0,\quad \text{and}\quad 0\le \frac{\partial}{\partial x_m}\Big\lvert_{0+}  h(x) < \infty,
\en
one has
\[
E^{Q_{u}}\left[ h\left(  X_{t_1}, \ldots, X_{t_m} \right) \right]=E^{Q_{s}}\left[ h\left(  X_{t_1}, \ldots, X_{t_m} \right) \right].
\] 

However, under both $Q^u$ and $Q^s$, the coordinate process is nonnegative almost surely. For nonnegative random variables the class of functions in \eqref{coffun} is a distribution determining class since indicators of open sets in $(\rr^{+})^m$ can be approximated by smooth functions satisfying \eqref{coffun}. This shows consistency and part (i) of the Proposition is proved. 

For part (ii) we use the change of measure in \eqref{abscont} to claim that
\[
E^R X_t = E^P 1/ X_t = E^Q \left(  X_t 1\{ X_t > 0 \} \frac{1}{X_t}  \right)= Q\left( X_t > 0 \right) = Q\left(  \tau_0 > t \right).
\]
If $Q(\tau_0 < \infty) > 0$, then $E^R X_t$ is not constant over $t$ which prevents it form being a martingale. 
This completes the proof of the Proposition.
\end{proof}

\section{Stochastic calculus with strict local martingales}

\comment{
Inspired by the previous representation theorem, we make the following definition:
\bigskip

\noindent\textbf{Definition.} We call an ordered pair probability
measures $(R,Q)$, defined on the canonical sample space of
continuous paths, a {Girsanov pair} if
\begin{enumerate}
\item under $R$, the coordinate process $X_t$ is a positive strict
local martingale starting from one;
\item under $Q$ the process $X_t$ is a nonnegative martingale;
\item The laws  $R$ and $Q$ are related by Proposition \ref{htrans}.
\end{enumerate}
}

The advantage of Proposition \ref{htrans} is that it allows us to
transport stochastic calculus with respect to strict local
martingales to that with actual martingales via a change of measure.

\begin{prop}\label{convexfn}
Let $R$ and $Q$ be two probability measures on $C[0,\infty)$ such that the coordinate process $X$ is a positive strict local martingale under $R$, is a nonnegative martingale under $Q$, and the relationship \eqref{abscont} holds.

Let $\tau_0$ denote the hitting time of zero of the coordinate process $X_t$.
Consider any nonnegative function $h:(0,\infty)\rightarrow \rr^+$.
For any bounded stopping time $\tau$, we get
\begin{equation}\label{convex}
E^R\left( h(X_{\tau}) \right)= E^Q g(X_{\tau})1_{\left\{ \tau_0 > \tau \right\}}.
\end{equation}
Here $g$ is the function $g(x)= xh({1}/{x})$, for all $x >0$.

Now suppose $\lim_{x\rightarrow 0}g(x)=\eta < \infty$. Define a map
$\gbar:[0,\infty)\rightarrow \rr$ by extending $g$ continuously,
i.e., $\gbar(x)=g(x)$ for $x >0$, and $\gbar(0)=\eta$. Then we have
\begin{equation}\label{gendecomp}
E^R\left( h(X_{\tau}) \right)= E^Q \gbar(X_{\tau}) - \eta Q\left( \tau_0 \le \tau\right).
\en
\end{prop}

\begin{proof}[Proof of Proposition \ref{convexfn}]
From \eqref{abscont} one gets
\[
E^R h(X_{\tau})= E^Q X_{\tau} h\left(\frac{1}{X_{\tau}}\right)
1_{\left\{ \tau_0 > \tau \right\}}= E^Q g(X_{\tau})1_{\left\{ \tau_0 > \tau \right\}},
\]
which proves \eqref{convex}.

For the second assertion note that for any nonnegative path $\omega$
which gets absorbed upon hitting zero, the following is an algebraic
identity:
\[
g(\omega_\tau)1_{\left\{ \tau_0 > \tau \right\}}= \gbar(\omega_{\tau}) -
\eta 1_{\left\{\tau_0 \le \tau  \right\}}.
\]
In particular this identity holds pathwise when $\omega$ is a path
of a nonnegative martingale. Thus we obtain
\[
\e^Q g(X_{\tau})1_{\left\{ \tau_0 > \tau \right\}}= \e^Q \gbar(X_{\tau}) - \eta Q(\tau_0 \le \tau).
\]
This proves \eqref{gendecomp}.
\end{proof}

As an example note that when
$h(x)=(x-a)^+$ for some $a \ge 0$, we get \eq\label{heaviside}
E^R\left( X_{\tau} - a \right)^+ = E^Q\left( 1 - aX_{\tau} \right)^+
- Q(\tau_0\le \tau). \en
We have the following corollary.

\begin{corollary}\label{corinfinity}
Let $h:(0,\infty) \rightarrow (0,\infty)$ be a function which is sublinear at infinity, i.e.,
\eq\label{sublin}
\lim_{x\rightarrow \infty}\frac{h(x)}{x}=0.
\en
Then, for all bounded stopping times $\tau$, one has
\eq\label{hgsublin}
E^R h(X_{\tau})= E^Q \gbar(X_{\tau}),\quad \gbar(x)=xh(1/x), \; x>0, \; \gbar(0)=0.
\en
\end{corollary}

\begin{proof}
The first part follows directly from \eqref{gendecomp} since
\[
\eta=\lim_{x\rightarrow 0}xh\left(\frac{1}{x}\right)=\lim_{x\rightarrow \infty}\frac{h(x)}{x}=0.
\]
The second conclusion is obvious.
\end{proof}

The previous corollary has some interesting consequences. For
example, when $h$ is convex, it is not difficult to verify that so
is $\gbar$. And thus both $h(X)$ and $\gbar(X)$ are submartingales
(under $R$ and $Q$ respectively) by \eqref{hgsublin}. This is in
spite of the strictness in the local martingale property of the
coordinate process under $R$. Additionally, if $h$ is symmetric with
respect to inverting $x$, i.e. $h(x)=xh(1/x)$, then $\gbar=h$. Hence
the strictness of local martingales has no effect when these functions
are applied. For example, $E^R \sqrt{X_{\tau}} = E^Q
\sqrt{X_{\tau}}$.

Before we end, let us mention that similar results can be
obtained from the semimartingale decomposition formulas of Madan and
Yor \cite{MY} and by PDE methods, as in Ekstr\"om and Tysk \cite{ET}.

\section{Convex functions of strict local martingales}

Strict local martingales are known for odd behavior which is not
shared by martingales. For example, a convex function of a
martingale is always a submartingale. This need not be the case with
local martingales.  However
if $N$ is a nonnegative strict local martingale, and $h$ is a convex
function sublinear at infinity, then $h(N_t)$ is again a
submartingale. This is in contrast to functions which are linear at
infinity. For example, in the case of $h(x)=x$, the process is
actually a supermartingale. Here we demonstrate another example, for
the function $(x-K)^+$ with $K >0$ applied to the inverse Bessel strict local martingale in Example 1. For a complementary result, see \cite{YY}.


\begin{prop}\label{besanal} Let $X_t$ be a BES($3$) process starting from one.
For any real $K \in [0, 1/2]$, the function $t\mapsto \e\{(1/X_t -
K)^+\}$ is strictly decreasing for all $t\in (0,\infty)$. However,
if $K > 1/2$, the function $t \mapsto \e\{(1/X_t - K)^+\}$ is
increasing in a neighborhood of zero but is strictly decreasing for all large enough $t$. 
In particular, it is strictly decreasing for
\[
t \ge \left( K\log\frac{2K+1}{2K-1} \right)^{-1}.
\]
\end{prop}


\begin{proof}[Proof of Proposition \ref{besanal}]

Let $B$ be a one-dimensional Brownian motion starting from one and absorbed
at zero. From the change of measure identity in \eqref{heaviside} we get
\[
\begin{split}
h(t)&:=\e\left\{\left( {1}/{X_t} - K\right)^+\right\} = \e (1-KB_{t\wedge \tau_0})^+ - P(\tau_0\le t).
\end{split}
\]
where $\tau_0$ is the hitting time of zero for the Brownian motion $B$.

Let $L_{t\wedge \tau_0}^{1/K}$ denote the local time at $1/K$ for the process $B_{\cdot \wedge \tau_0}$.
If we take derivatives with respect to $t$ in the equation above, we get
\eq\label{firstderiv}
\begin{split}
h'(t)&= \frac{d}{dt}E(1-KB_{t\wedge \tau_0})^+ - \frac{d}{dt}P(\tau_0\le t)\\
&= \frac{K}{2}\frac{d}{dt}E L_{t\wedge \tau_0}^{1/K} - \frac{d}{dt}P(\tau_0\le t).
\end{split}
\en
The second equality above is due to the Tanaka formula.

Now to compute the first term on the right of \eqref{firstderiv}, we use \cite[Ex.~1.12, p.~407]{ry99} and the explicit (see \cite[page 97]{KS}) transition function of Brownian motion absorbed at zero to get
\[
\begin{split}
\frac{d}{dt}E L_{t\wedge \tau_0}^{1/K}=\frac{1}{\sqrt{2\pi t}}\left[ e^{-(1-1/K)^2/2t} -
e^{-(1+1/K)^2/2t}\right] .
\end{split}
\]

The second term on the right side of \eqref{firstderiv} above is
the density of the first hitting time of zero, which we know
(\cite[page 80]{KS}) to be $({2\pi t^3})^{-1/2}e^{-1/2t}$.

\comment{
\begin{lem}\label{mkovderiv} Suppose $\{X_t, \; t \ge 0\}$ is a continuous nonnegative local
martingale which satisfies the following SDE
\begin{equation}\label{sdeform}
d X_t= \sigma(t,X_t)d\beta_t, \quad t\in [0,\infty),\quad X_0=1.
\end{equation}
Here $\beta$ is a one-dimensional standard Brownian motion and
$\sigma(t,x)$ is some measurable nonnegative function on
$\rr^+\times \rr^+$.

Further assume that the process $X_t$ admits a continuous marginal
density at each time $t$ at every strict positive point $y$ which is
given by
\[
p_t(y)= P\left( X_t\in dy\;\Big\lvert\; X_0= 1\right), \quad y > 0.
\]
Let $L_t^a$ denote the local time of $X$ at level $a > 0$ and at time
$t$. Then
\begin{equation}\label{derivlocal}
\frac{d}{dt}\e\left(L_t^a\right)= \sigma^2(t,a)p_t(a).
\end{equation}
\end{lem}
}
 
 Thus, we get
\eq\label{deriv2}
h'(t)= \frac{K}{2\sqrt{2\pi t}}\left[ e^{-(1-1/K)^2/2t} -
e^{-(1+1/K)^2/2t}\right] - \frac{1}{\sqrt{2\pi t^3}} e^{-1/2t}.
\en
Thus, $h'(t) < 0$ if and only of
\eq\label{deriv3}
\begin{split}
\frac{2}{K t} &> e^{1/2t}\left[ e^{-(1-1/K)^2/2t} - e^{-(1+1/K)^2/2t}\right]\\
&=\exp\left[\frac{(2K-1)}{2K^2t}\right] - \exp\left[-\frac{(2K+1)}{2K^2t} \right].
\end{split}
\en
If $K > 1/2$, it is clear that the above inequality cannot be true in a neighborhood of zero. Thus, the function $h(t)$ must be initially increasing. 

To get an estimate of the point when it starts decreasing, let $t=1/y$. Consider the function on
the right side of the last inequality. We consider two
separate cases. First suppose $K > 1/2$. Then both $2K-1$ and $2K+1$
are positive. If for two positive parameters $\lambda_2 > \lambda_1
>0$, we define a function $q$ by $q(y) = \exp(\lambda_1 y) -
\exp(-\lambda_2 y)$, $y >0$, it then follows that
\begin{equation}\label{whatisq}
\begin{split}
q'(y)&= \lambda_1e^{\lambda_1 y} + \lambda_2e^{-\lambda_2y},\quad q'(0)=\lambda_1+\lambda_2,\\
q^{\prime\prime}(y)&= \lambda_1^2e^{\lambda_1 y} -
\lambda_2^2e^{-\lambda_2y}.
\end{split}
\end{equation}
Note that $q^{\prime\prime}(y) < 0$, for all
\begin{equation}\label{whenneg}
0 \le y  < \frac{2\log\left(\lambda_2/\lambda_1\right)}{\lambda_1+\lambda_2}.
\end{equation}
Since $q'(y)$ is always positive, it follows that $q$ is an
increasing concave function starting from zero in the interval given
by \eqref{whenneg}. Thus,
\begin{equation}\label{compare2}
q(y)= q(y)-q(0) < yq'(0)= y\left( \lambda_1 + \lambda_2\right).
\end{equation}

Take $\lambda_1=(2K-1)/2K^2$ and $\lambda_2=(2K+1)/2K^2$. Then $\lambda_1 +
\lambda_2=2/K$. By \eqref{whenneg} we get that if
\[
y \le C_1:={K\log\frac{2K+1}{2K-1}},
\]
then, from \eqref{compare2} it follows
\[
K\left\{\exp\left[\frac{(2K-1)y}{2K^2}\right] - \exp\left[-\frac{(2K+1)y}{2K^2}
\right]\right\} < 2y.
\]
That is, by \eqref{deriv3}, $h'(t) < 0$, i.e., $h$ is strictly decreasing for all
\[
t > \left( K\log\frac{2K+1}{2K-1} \right)^{-1}.
\]

The case when $0 <K \le 1/2$ can handled similarly. Suppose $0 < \lambda_1 <
\lambda_2$ are positive constants. Consider
the function
\[
r(y)= -\lambda_1 y + \lambda_2 y - e^{-\lambda_1 y} + e^{-\lambda_2 y}, \quad y \in [0,\infty).
\]
Then $r(0)=0$, and
\[
r'(y)= -\lambda_1\left( 1- e^{-\lambda_1 y} \right) + \lambda_2\left(
1-e^{-\lambda_2 y}\right) > 0,
\quad y\in[0,\infty),
\]
because $\lambda_1 < \lambda_2$. Thus, for all positive $y$, we have $r(y) >0$, i.e.,
\[
e^{-\lambda_1 y} - e^{-\lambda_2 y} < (-\lambda_1 + \lambda_2) y.
\]
We use this for $\lambda_1= (1-2K)/2K^2$ and $\lambda_2=(1+2K)/2K^2$. Note that,
as before $(-\lambda_1 + \lambda_2)y= 2y/K$.

From \eqref{deriv3} it follows that $h'(t) > 0$ for all
$t\in(0,\infty)$. Thus we have established that if $K \le 1/2$, the
function $t\mapsto \e({1}/{X_t} - K)^+$ is strictly decreasing for
all $t\in(0,\infty)$. This completes the proof of the proposition.
\end{proof}

For mathematical completeness we show below that a similar result
can be proved for the Bessel process starting from zero, although in
this case there is no dependence on $K$. The proof is much simpler
and essentially follows by a scaling argument. Note that, even in
this case the reciprocal of the Bessel process is well-defined for
all times except at time zero. Hence $1/X_t$, $t\in (0,\infty)$, can
be thought as a Markov process with an \textit{entrance
distribution}, i.e., a pair consisting of a time-homogenous Markov
transition kernel $\{ P_{t} \}$, $t > 0$, and a family of
probability measures $\{ \mu_s \}$, $s > 0$, satisfying the
constraint $\mu_s\ast P_{t}= P_{t+s}$. Here $\ast$ refers to the
action of the kernel on the measure.

\begin{prop} Let $X_t$ be a $3$-dimensional Bessel process, BES(3), such that
$X_0=0$. For any two time points $u > t > 0$, and for $K \ge 0$, one
has
\begin{equation}\label{main}
\e\left(\frac{1}{X_u} - K\right)^+ < \e\left(\frac{1}{X_t} - K\right)^+.
\end{equation}
\end{prop}

\begin{proof}
Fix $u > t$. Recall that BES($3$), being the norm of a three
dimensional Brownian motion, has the Brownian scaling property when
starting from zero. That is to say, for any $c > 0$,
\[
\left(\;\frac{1}{\sqrt{c}}X_{cs},\; s \ge 0\;\right) \stackrel{\mcal{L}}{=}
\left(\;X_{s},\; s \ge 0\;\right),
\]
where the above equality is equality in law.

Take $c=u/t$, and apply the above equality for $X_s$ when $s=t$, to
infer that $c^{-1/2}X_u$ has the same law as $X_t$, and thus
\begin{equation}\label{main2}
\e\left(\frac{1}{X_u} - K\right)^+=\e\left(\frac{c^{-1/2}}{X_t} - K\right)^+=
c^{-1/2}\e\left(\frac{1}{X_t} - \sqrt{c}K\right)^+.
\end{equation}

Note that for any $\sigma >1$, we have $\left( x - \sigma
K\right)^+/\sigma < (x - K)^+,\quad \forall\; x > 0$. Since $c > 1$,
taking $\sigma=\sqrt{c}$, one deduces from \eqref{main2}
\[
\e\left(\frac{1}{X_u} - K\right)^+ < \e\left(\frac{1}{X_t} - K\right)^+,
\]
which proves the result.
\end{proof}

We conclude this subsection with an example of a strict local
martingale $S$ where $E(S_t - K)^+$ is not asymptotically decreasing
for any $K$. This, coupled with the earlier Bessel result,
establishes the fact that functions which are not sublinear at
infinity can display a variety of characteristics when applied to
strict local martingales.

We inductively construct a process in successive intervals
$[i,{i+1})$ by the following recipe. The process starts at zero. The
process in the odd interval $[2i, 2i+1)$ is an exponential Brownian
motion $\exp(B_t - t/2)$ starting from $S_{2i}$ and independent of
the past. On the even intervals $[2i+1,2i+2)$ the process $S$ is an
inverse Bessel process starting from $S_{2i+1}$ and again
independent of the past. The constructed process is always a
positive local martingale. The value of the function $E(S_t-K)^+$ is
increasing in the odd intervals due to the martingale component, and
decreasing (at least when $K\le 1/2$) in the Bessel component by
Proposition \ref{besanal}.

One might object to the fact that this process is not strict local
throughout. But, one can mix the two components, by a sequence of
coin tosses which decides whether to use Brownian or the Bessel
component in the corresponding interval. By choosing the probability
of heads in these coins in a suitably predictable manner, we can
generate a local martingale which is strict throughout but
$E(S_t-K)^+$ does not decrease anywhere.

\subsection{A multidimensional analogue by Kelvin transform}

In the last section we saw that for any strict local martingale law
$R$ there is a true martingale law $Q$ such that equality
\eqref{convex} holds. The transformation $g(x)=xh(1/x)$ has a
well-known analogue in dimensions higher than two called the Kelvin
transform. In our final subsection we present an interesting
multidimensional generalization of our results for dimensions $d >
2$. We take the following definition from the excellent book on
harmonic function theory \cite[Chapter 4]{harmonicfn} by Axler,
Bourdon, and Ramey.

The Kelvin transform $K$ is an operator acting on the space of real
functions $u$ on a subset of $\rr^d\backslash \{0\}$. Let $u$ be a
$C^2$ function on an open subset $\mathcal{D}$ of $\rr^d\backslash
\{0\}$. Let $\mathcal{D}^*$ be the image of $\mathcal{D}$ under the
inversion map
\eq\label{invmap}
x\mapsto x^*=x/\abs{x}^2.
\en
For such a $u$, we define a function $K[u]:\mathcal{D}^* \rightarrow \rr$ by the formula
\eq\label{whatisK}
K[u](y)=\abs{y}^{2-d}u\left( {y}/{\abs{y}^2} \right).
\en
Notice that $K$ is its own inverse.

The most striking property of this transform is that $K$
\textit{commutes} with the Laplacian (\cite[page 62]{harmonicfn}).
Let $v$ be the function $v(x)=\abs{x}^4 \Delta u(x), \; x \in
\mathcal{D}$. Then, at any point $y\in \mathcal{D}^*$, we have
\eq\label{kelvincomm} \Delta K[u](y)= K\left[v\right](y). \en In
particular, if $u$ is harmonic in $\mathcal{D}$ (i.e., $\Delta
u=0$), then $K[u]$ is harmonic in $\mathcal{D}^*$. Also, if $u$ is
subharmonic (i.e., $\Delta u \ge 0$), then so is $K[u]$.

A $d$-dimensional conformal local martingale is a
process $(X(1), \ldots, X(d))$ such that each coordinate $X(i)$ is a
local martingale and
\[
\langle X(i), X(j) \rangle = \iprod{X(1)}1_{\{i=j\}}, \quad \text{for
all}\; 1\le i,j \le d.
\]
We consider it to be \textit{strict local} if at least one of its
coordinate processes is a strict local martingale. We replace the
condition of nonnegativity of one dimensional processes by
restricting our multidimensional processes in the complement of a
compact neighborhood of zero.

\begin{prop}\label{multiprop}
Let $D$ be a compact neighborhood of the origin and denote its
complement $\rr^d\backslash D$ by $D^c$. Let $P$ be a probability
measure on the canonical sample space $C^d[0,\infty)$ such that, under $P$, the
coordinate process $X$ is a conformal local martingale in $\rr^d$
($X_0=x_0\in D^c$) absorbed upon hitting the boundary of $D$. Then
there is a probability measure $Q$ such that, under $Q$, the
coordinate process is a true conformal martingale which takes values
in $(D^c)^*$ (see \eqref{invmap}) such that, for any bounded measurable function $U:D^c \rightarrow \rr$, and
for any bounded stopping time $\tau$, we have \eq\label{changemulti}
\abs{x_0}^{2-d}E^P\left[ U\left( X_{\tau} \right)  \right]=
E^Q\left[  \abs{{X_\tau}}^{2-d} U\left(
\frac{X_{\tau}}{\abs{X_{\tau}}^2} \right) \right]. \en
\end{prop}

\comment{
\begin{prop}\label{mainpropmulti}
Let $D$ be a compact neighborhood of the origin and denote its
complement $\rr^d\backslash D$ by $D^c$. Let $X$ be any conformal
local martingale in $\rr^d$ ($X_0\in D^c$) absorbed upon hitting the
boundary of $D$. Let $U$ be a $C^2$-subharmonic function on the
complement of the open disc, i.e., $\Delta U \ge 0$ on $\{x:\;
\abs{x}\ge 1\}$. Assume more over that $U$ is  {harmonic at
infinity} in the sense that $\lim_{x\rightarrow \infty}U(x)=0$. Then
$U(X)$ is a submartingale.
\end{prop}
}

The proof requires the following lemmas.

\begin{lem}
Let $\abs{\cdot}$ denote the Euclidean norm in dimension $d$. Let
$\tau_1$ be the hitting time of $D$ i.e.,
\[
\quad \tau_1=\inf\left\{t\ge 0:\; X_t\in D\right\}.
\]
Then, the process $\abs{X_{t\wedge \tau_1}}^{2-d}, \; t\ge 0$, is a $P$-martingale.
\end{lem}

\begin{proof}
The function $\abs{x}^{2-d}$ is harmonic in $\rr^d$. Thus
$\abs{X_t}^{2-d}$ is a local martingale itself. Since it is bounded
in $D^c$, it must be a true martingale.
\end{proof}

Suppose $X_0=x_0\in D^c$. We can change the law of $X$ by using
$\abs{X_{t\wedge \tau_1}}^{2-d}$ as a Radon-Nikod\'ym derivative
(after normalizing). We get the following lemma.

\begin{lem}\label{inversionmgle}
Suppose $X_0=x_0$ such that $x_0\in D^c$. We change $P$ by using the positive martingale
\[
\phi(X_t)=\abs{X_{t\wedge \tau_1}}^{2-d}/\abs{x_0}^{2-d}
\]
as a Radon-Nikod\'ym derivative. Call this measure $Q^*$. Under
$Q^*$, the process \eq\label{whatisy} Y_t=\frac{X_{t\wedge
\tau_1}}{\abs{X_{t\wedge \tau_1}}^2} \en is again a $d$-dimensional
conformal local martingale such that every coordinate process $Y$ is
a true martingale.
\end{lem}

\begin{proof}[Proof of Lemma \ref{inversionmgle}]
Note that, since $D$ is a neighborhood of zero, the set $(D^c)^*$ is
compact. Thus, $Y$ lives in a bounded set. To show that the $i$th
coordinate process $Y(i)$ is a local martingale, we use the harmonic
function
\[
u(x)=x(i)\quad \text{on}\quad \Omega=\rr^n \backslash \{0\}.
\]
By \eqref{kelvincomm}, its Kelvin transform is also harmonic. Hence the process
\[
\abs{X_t}^{2-d}u\left( X_t/\abs{X_t}^2 \right)= \abs{X_t}^{-d}X_t(i)
\]
is a local martingale under $P$. But, since $\abs{X_t}^{2-d}$ (until
$\tau_1$) is used as a Radon-Nikod\'ym derivative after being
scaled, by Bayes rule \cite[page 193]{KS}, the process $u\left(
X_{t\wedge \tau_1}/\abs{X_{t\wedge \tau_1}}^2 \right)$ is a local
martingale under the changed measure $Q^*$. But this implies that
$Y(i)$ is a local martingale under $Q^*$. But, since every $Y(i)$ is
bounded they must be true martingales.

To show that $Y$ is conformal, we use the harmonic function
$u(x)=x(i)x(j)$ again on the full domain $\rr^n\backslash \{0\}$ for
any pair of coordinates $i\neq j$. Exactly as in the previous
paragraph, we infer from \eqref{kelvincomm} that $Y_t(i)Y_t(j)$ is a
local martingale under $Q$. But this implies
$\iprod{Y(i),Y(j)}\equiv 0$. That their quadratic variations must be
the same follows from symmetry.
\end{proof}

\begin{proof}[Proof of Proposition \ref{multiprop}]
The construction of $Q$ has been done in Lemma \ref{inversionmgle}
where it is the law of the process $Y$ under $Q^*$. Note that the
Radon-Nikod\'ym derivative $\phi(X_t)$ never hits zero under $P$. Thus, even
under $Q^*$, the process $Y$ never hits zero. To show the equality
\eqref{changemulti}, we use the change of measure to get
\[
\begin{split}
E^{Q^*}\left[  \abs{Y_{\tau}}^{2-d} U\left(
\frac{Y_{\tau}}{\abs{Y_{\tau}}^2} \right)   \right]&=\abs{x_0}^{d-2}E^P\left[ \abs{X_{\tau}}^{2-d}\abs{X_{\tau}}^{d-2}
U\left(  X_{\tau} \right) \right]\\
&=\abs{x_0}^{d-2}E^P\left[ U\left(  X_{\tau} \right) \right].
\end{split}
\]
This completes the proof of the result.
\end{proof}

\section{Applications to financial bubbles}\label{discussfin}

A natural question is: what happens to a financial market when the
no arbitrage condition yields a strict local martingale (rather than
a true martingale) under a risk neutral measure? Several authors
have looked at this problem and offered solutions to anomalies which
might result from the lack of the martingale property. One
interesting perspective offered in this direction is the theory of
price bubbles as argued in 2000 by Loewenstein and Willard
~\cite{HLW}. They propose that to identify a bubble one needs to
look at the difference between the market price of an asset and its
fundamental price. Their argument is later complemented and further
developed by Cox and Hobson \cite{CH} and the two articles by
Jarrow, Protter, and Shimbo \cite{JPS1}, \cite{JPS2}. 
\comment{Please see the
latter articles for the definitions of the market and the
fundamental prices of an asset and any of the other financial terms
that follow. In particular, the authors in \cite{JPS1} and
\cite{JPS2} classify bubbles into three types in an arbitrage-free
market satisfying Merton's \textit{No Dominance} condition (see
\cite{JPS1} or \cite{merton73}). One, in which the difference
between the two price processes under an equivalent local martingale
measure is a uniformly integrable martingale; two, when it is a
martingale but non-uniformly integrable; and last, when it is a
strict local martingale. In a static market with infinite horizon,
for a stock which pays no dividends, Example 5.4 in \cite{JPS2}
shows that the difference between the two prices is actually the
current market price of the stock.}
In particular they identify that a stock price which behaves
as a strict local martingale under an equivalent local martingale
measure is an example of a price bubble. Moreover, Cox and
Hobson \cite{CH} go on to exhibit (among
other things) how in the presence of bubbles put-call parity might
not hold and call prices do not tend to zero as strike tends to
infinity.

We consider a market with a single risky asset (stock) and zero spot
interest rate. Let $\{S_t\}$, $t\in(0,\infty)$, be a positive
continuous strict local martingale which models the discounted price
of the (non-dividend paying) stock under an equivalent local
martingale measure. We have the following result which follows
immediately from Corollary \ref{corinfinity} and the subsequent
Bessel example.

\begin{prop}
Suppose for a European option, the discounted pay-off at time $T$ is
given by a convex function $h(S_T)$ which is sublinear at infinity,
i.e., $\lim_{x\rightarrow \infty}h(x)/x=0$. Then the price of the
option is increasing with the time to maturity, $T$, whether or not
a bubble is present in the market. In other words, $E(h(S_T))$ is an
increasing function of $T$. For example, consider the put option
with a pay-off $(K-x)^+$.

However, for a European call option, the price of the option
$E(S_T-K)^+$ with strike $K$ might decrease as the maturity
increases.
\end{prop}

This feature may seem strange at first glance, but if we assume the
existence of a financial bubble, the intuition is that it is
advantageous to purchase a call with a short expiration time, since
at the beginning of a bubble prices rise, sometimes dramatically.
However in the long run it is disadvantageous to have a call,
increasingly so as time increases, since the likelihood of a crash
in the bubble taking place increases with time. 

It is worthwhile to compare our results with processes that are increasing in the convex order. 
We refer the interested reader to the articles by Hirsch and Yor \cite{HY091}, \cite{HY092} and the references within.

Of course, pricing a European option by the usual formula when the
underlying asset price is a strict local martingale is itself
controversial. For example, Heston, Loewenstein, and Willard
\cite{helowill06} observe that under the existence of bubbles in the
underlying price process, put-call parity might not hold, American
calls have no optimal exercise policy, and lookback calls have
infinite value. Madan and Yor~\cite{MY} have recently argued that
when the underlying price process is a strict local martingale, the
price of a European call option with strike rate $K$ should be
modified as $\lim_{n\rightarrow \infty}E\left[(S_{T\wedge T_n} -
K)^+ \right]$, where $T_n= \inf\left\{ t\ge 0: S_t \ge n \right\}$,
$n \in \mathbb{N}$, is a sequence of hitting times. This proposal
does however, in effect, try to hide the presence of a bubble and
act as if the price process is a true martingale under the risk
neutral measure, rather than a strict local martingale.

A different approach to such market
anomalies has been studied extensively by several authors
when the
\emph{candidate Radon-Nikod\'ym derivative} for the risk-neutral
measure turns out to be a strict local martingale.
This is the spirit taken in Platen \cite{P}, Platen and Heath \cite{PH}, Fernholz and Karatzas
\cite{ferkar05}, Fernholz, Karatzas, and Kardaras
\cite{ferkarkar05}, and Karatzas and Kardaras \cite{karkar07}. For the latter group of authors this is intimately connected
with the concept of a \textit{weakly diverse market} and displays a number of anomalies similar to the case of bubbles. For
example, put-call parity fails to hold in such markets. See, Remark
9.1 and 9.3, and Example 9.2 in \cite{ferkarkar05}.  Finally, we mention that P. Guasoni and M. R\'asonyi have recently studied the robustness of bubble models by showing that if $S$ is a price process which models a financial bubble, under some conditions one can find another process $S^\epsilon$ which is ``$\epsilon$ close to $S$'' and is not itself a bubble.  Therefore in this sense models for bubbles perhaps can be seen as delicate from the standpoint of model risk~\cite{GR}.

\section*{Acknowledgements.} We are grateful to Professors
Marc Yor and Monique Jeanblanc who drew our attention to the papers
of Madan and Yor~\cite{MY} and Elworthy, Li and Yor~\cite{ELY}.  We
also thank the anonymous referees for their excellent reports on a
previous version of this paper. The second author gratefully
acknowledges benefitting from a Fulbright-Tocqueville Distinguished
Chair award at the University of Paris -- Dauphine, during the
development of this research.

\end{document}